\documentclass[12pt,a4paper]{article}
\usepackage[cp1251]{inputenc} 
\usepackage{amsfonts}

\newcommand{\di}{\displaystyle}
\newcommand{\B}{$\hfill\Box$}
\newcommand{\al}{\alpha}

\newcommand{\ga}{\gamma}
\newcommand{\de}{\delta}
\newcommand{\la}{\lambda}
\newcommand{\om}{\omega}
\newcommand{\ee}{\varepsilon}
\newcommand{\vv}{\varphi}
\newcommand{\iy}{\infty}

\begin{document}

\begin{center}
{\large\bf
Recovering Higher Order Differential Systems with Regular Singularities on Star-type Graphs}\\[0.2cm]
{\bf V.\,Yurko} \\[0.2cm]
\end{center}

\thispagestyle{empty}

{\bf Abstract.} We study an inverse spectral problem for arbitrary order ordinary differential equations on compact star-type graphs when differential equations have regular singularities at boundary vertices. As the main spectral characteristics we introduce and study the so-called Weyl-type matrices which are generalizations of the Weyl function (m-function) for the classical Sturm-Liouville operator. We provide a procedure for constructing the solution of the inverse problem and prove its uniqueness.



{\bf 1. Introduction. } We study an inverse spectral problem for
arbitrary order ordinary differential equations on compact
star-type graphs when differential equations have regular
singularities at boundary vertices. Boundary value problems on
graphs (spatial networks, trees) often appear in natural sciences
and engineering (see [1-6]). Inverse spectral problems consist in
recovering operators from their spectral characteristics. We pay
attention to the most important nonlinear inverse problems of
recovering coefficients of differential equations (potentials)
provided that the structure of the graph is known a priori.

For {\it second-order} differential operators on compact graphs
inverse spectral problems have been studied fairly completely in
[7-13] and other works. Inverse problems for {\it higher-order}
differential operators on graphs were investigated in [14-15]. We
note that inverse spectral problems for second-order and for
higher-order ordinary differential operators on {\it an interval}
have been studied by many authors (see the monographs [16-22] and
the references therein). Arbitrary order differential operators on
an interval with regular singularities were considered in [23-26].

In this paper we study the inverse spectral problem for arbitrary
order differential operators with regular singularities on compact
star-type graphs. As the main spectral characteristics in this
paper we introduce and study the so-called Weyl-type matrices
which are generalizations of the Weyl function (m-function) for
the classical Sturm-Liouville operator (see [27]), of the Weyl
matrix for higher-order differential operators on an interval
introduced in [21-22], and generalizations of the Weyl-type
matrices for higher-order differential operators on graphs (see
[14-15]). We show that the specification of the Weyl-type matrices
uniquely determines the coefficients of the differential equation
on the graph, and we provide a constructive procedure for the
solution of the inverse problem from the given Weyl-type matrices.
For studying this inverse problem we develope  the method of
spectral mappings [21-22]. We also essentially use ideas from [23]
on differential equations with regular singularities. The obtained
results are natural generalizations of the well-known results on
inverse problems for differential operators on an interval and on graphs.\\

{\bf 2. Formulation of the inverse problem. } Consider a compact
star-type graph $T$ in ${\bf R^\om}$ with the set of vertices
$V=\{v_0,\ldots, v_p\}$ and the set of edges ${\cal
E}=\{e_1,\ldots, e_p\},$ where $v_1,\ldots, v_{p}$ are the
boundary vertices, $v_0$ is the internal vertex, and
$e_j=[v_{j},v_0],$ $j=\overline{1,p},$
$\di\bigcap_{j=1}^{p}\,e_j=\{v_0\}$. Let $l_j$ be the length of
the edge $e_j$. Each edge $e_j\in {\cal E}$ is parameterized by
the parameter $x_j\in [0,l_j]$ such that $x_j=0$ corresponds to
the boundary vertices $v_1,\ldots, v_{p}$, and $x_j=l_j$
corresponds to the internal vertex $v_0$. An integrable function
$Y$ on $T$ may be represented as $Y=\{y_j\}_{j=\overline{1,p}}$,
where the function $y_j(x_j)$ is defined on the edge $e_j$.

Consider the differential equations on $T$:
$$
y_j^{(n)}(x_j)+\di\sum_{\mu=0}^{n-2}
\Big(\frac{\nu_{\mu j}}{x_j^{n-\mu}}+ q_{\mu j}(x_j)\Big)y_j^{(\mu)}(x_j)
=\la y_j(x_j),\quad x_j\in (0, l_j),\quad j=\overline{1,p},                          \eqno(1)
$$
where $\la$ is the spectral parameter, $q_{\mu j}(x_j)$ are complex-valued integrable
functions. We call $q_j=\{q_{\mu j}\}_{\mu=\overline{0,n-2}}$ the potential on the edge
$e_j$, and we call $q=\{q_{j}\}_{j=\overline{1,p}}$ the potential on the graph $T.$
Let $\{\xi_{kj}\}_{k=\overline{1,n}}$ be the roots of the characteristic polynomial
$$
\de_j(\xi)=\sum_{\mu=0}^n \nu_{\mu j} \prod_{k=0}^{\mu-1} (\xi-k),
\quad \nu_{nj}:=1,\; \nu_{n-1,j}:=0.
$$
For definiteness, we assume that $\xi_{kj}-\xi_{mj}\ne sn, s\in{\bf Z},$
$Re\,\xi_{1j}<\ldots <Re\,\xi_{nj},$ $\xi_{kj}\ne\overline{0,n-3}$ (other cases
require minor modifications). We set $\theta_j:=n-1-Re\,(\xi_{nj}-\xi_{1j}),$
and assume that the functions $q_{\mu j}^{(\nu)}(x_j),$ $\nu=\overline{0,\mu-1},$
are absolutely continuous, and $q_{\mu j}^{(\mu)}(x_j)x_j^{\theta_j}\in L(0,l_j).$

\medskip
Let $\la=\rho^{n},$ $\ee_k=\exp(2\pi ik/n),$ $k=\overline{0,n-1}.$
It is known that the $\rho$ -- plane can be partitioned into sectors $S$ of angle
$\frac{\pi}{n}$ $\Big(\arg\rho\in\Big(\frac{k_0\pi}{n},\frac{(k_0+1)\pi}{n}\Big),$
$k_0=\overline{-n,n-1}\Big)$ in which the roots $R_{1}, R_{2},\ldots, R_{n}$ of
the equation $R^{n}-1=0$ can be numbered in such a way that
$$
Re(\rho R_{1})<Re(\rho R_{2})<\ldots<Re(\rho R_{n}),\quad\rho\in S.                \eqno(2)
$$
Clearly, $R_k=\ee_{\eta_k}$, where $\eta_1,\ldots, \eta_n$ is a permutation of
the numbers $0,1,\ldots, n-1,$ depending on the sector. Let us agree that
$$
\rho^{\mu}=\exp(\mu(\ln|\rho|+i\arg\,\rho)),\; \arg\,\rho\in(-\pi,\pi],
\quad R_k^{\mu}=\exp(2\pi i\mu\eta_k/n).
$$
Let the numbers $c_{kj0},\; k=\overline{1,n},$ be such that
$$
\prod_{k=1}^n c_{kj0}=
\Big(\det[\xi_{kj}^{\nu-1}]_{k,\nu=\overline{1,n}}\Big)^{-1}.
$$
Then the functions
$$
C_{kj}(x_j,\la)=x_j^{\xi_{kj}} \sum_{\mu=0}^\iy c_{kj\mu} (\rho x_j)^{n\mu},
\quad c_{kj\mu}=c_{kj0}\Big(\prod_{s=1}^\mu \de_j(\xi_{kj}+sn)\Big)^{-1},
$$
are solutions of the differential equation in the case when
$q_{\mu j}(x_j)\equiv 0,$ $\mu=\overline{0,n-2}.$ Moreover,
$\det[C_{kj}^{(\nu-1)}(x_j,\la)]_{k,\nu=\overline{1,n}}\equiv 1.$
Denote $\rho^*=2n\max\limits_{\mu,j} \|q_{\mu j}\|_{L(0,l_j)},\;
\mu=\overline{0,n-2},\; j=\overline{1,p}.$ In [23] we constructed
 special fundamental systems of solutions
$\{S_{kj}(x_j,\la)\}_{k=\overline{1,n}}$ and
$\{E_{kj}(x_j,\rho)\}_{k=\overline{1,n}}$ of equation (1) on the
edge $e_j$, possessing the following properties.

1) For each $x_j\in(0,l_j],$ the functions
$S_{kj}^{(\nu)}(x_j,\la),\;\nu=\overline{0,n-1},$ are entire in
$\la.$ For each fixed $\la,$ and $x_j\to 0,$
$$
S_{kj}(x_j,\la)\sim c_{kj0} x_j^{\xi_{kj}},\quad
(S_{kj}(x_j,\la)-C_{kj}(x_j,\la))x_j^{-\xi_{kj}}=o(x_j^{\xi_{nj}-\xi_{1j}}).
$$
Moreover, $\det[S_{kj}^{(\nu-1)}(x_j,\la)]_{k,\nu=\overline{1,n}}\equiv 1,$ and
$|S_{kj}^{(\nu)}(x_j,\la)|\le C|x^{\xi_{kj}-\nu}|,\; |\rho|x_j\le 1.$
Here and below, we shall denote by the same symbol $C$ various positive constants
in the estimates independent of $\la$ and $x_j$.

2) For each $x_j>0$ and for each sector $S$ with property (2), the functions
$E_{kj}^{(\nu)}(x_j,\rho),$ $\nu=\overline{0,n-1},$ are regular with respect to
$\rho\in S,\; |\rho|>\rho^*$, and continuous for $\rho\in \overline{S},\;
|\rho|\ge \rho^*$. Moreover,
$$
|E_{kj}^{(\nu)}(x_j,\rho)(\rho R_k)^{-\nu}\exp(-\rho R_k x_j)-1|
\le C(|\rho|x_j),\quad \rho\in \overline{S},\quad |\rho|x_j\ge 1.
$$

3) The relation
$$
E_{kj}(x_j,\rho)=\sum_{\mu=1}^n b_{kj\mu}(\rho)S_{\mu j}(x_j,\la),           \eqno(3)
$$
holds, where
$$
b_{kj\mu}(\rho)= b_{\mu j}^0 R_k^{\xi_{\mu j}}\rho^{\xi_{\mu
j}}[1],\quad b_{\mu j}^0\ne 0,\quad \rho\in\overline{S},\quad
\rho\to\iy,            \eqno(4)
$$
$$
\prod_{\mu=1}^n b_{\mu j}^0=\det[R_k^{\nu-1}]_{k,\nu=\overline{1,n}}
\Big(\det[R_k^{\xi_{\mu j}}]_{k,\mu=\overline{1,n}}\Big)^{-1},
$$
where $[1]=1+O(\rho^{-1}).$

Note that the asymptotical formula (4) is the most important and
nontrivial property of these solutions. This property allows one
to study both direct and inverse problems for arbitrary order
differential operators with regular singularities (see [24-26]).

Consider the linear forms
$$
U_{j\nu}(y_j)=\sum_{\mu=0}^{\nu}\ga_{j\nu\mu}y_j^{(\mu)}(l_j),
\; j=\overline{1,p},\; \nu=\overline{0,n-1},
$$
where $\ga_{j\nu\mu}$ are complex numbers, $\ga_{j\nu}:=\ga_{j\nu\nu}\ne 0.$
The linear forms $U_{j\nu}$ will be used in matching conditions at the
internal vertex $v_0$ for boundary value problems and for the correspondung
special solutions of equation (1).

\medskip
Fix $s=\overline{1,p},\; k=\overline{1,n-1}.$ Let $\Psi_{sk}=
\{\psi_{skj}\}_{j=\overline{1,p}}$ be solutions of equation (1)
on the graph $T$ under the boundary conditions
$$
\psi_{sks}(x_s,\la)\sim c_{ks0}x_s^{\xi_{ks}},\quad x_s\to 0,             \eqno(5)
$$
$$
\psi_{skj}(x_j,\la)=O(x_j^{\xi_{n-k+1,j}}),
\quad x_j\to 0,\quad j=\overline{1,p},\;j\ne s,                           \eqno(6)
$$
and the matching conditions at the vertex $v_0$:
$$
U_{1\nu}(\psi_{sk1})=U_{j\nu}(\psi_{skj}),
\;\; j=\overline{2,p},\; \nu=\overline{0,k-1},                           \eqno(7)
$$
$$
\di\sum_{j=1}^{p} U_{j\nu}(\psi_{skj})=0,\quad \nu=\overline{k,n-1}.     \eqno(8)
$$
The function $\Psi_{sk}$ is called the Weyl-type solution of order $k$
with respect to the boundary vertex $v_s$. Define additionally
$\psi_{sns}(x_s,\la):=S_{ns}(x_s,\la).$

Using the fundamental system of solutions $\{S_{\mu j}(x_j,\la)\}$
on the edge $e_j$, one can write
$$
\psi_{skj}(x_j,\la)=\di\sum_{\mu=1}^{n}
M_{skj\mu}(\la)S_{\mu j}(x_j,\la),
\quad j=\overline{1,p}, \quad k=\overline{1,n-1},                       \eqno(9)
$$
where the coefficients $M_{skj\mu}(\la)$ do not depend on $x_j.$

It follows from (9) and the boundary condition (6) for the
Weyl-type solutions that
$$
\psi_{skj}(x_j,\la)=\sum_{\mu=n-k+1}^{n} M_{skj\mu}(\la)S_{\mu
j}(x_j,\la), \quad j=\overline{1,p}\setminus s. \eqno(10)
$$
Similarly, using (5) one gets
$$
\psi_{sks}(x_s,\la)=S_{ks}(x_s,\la)+\sum_{\mu=k+1}^{n}
M_{sk\mu}(\la)S_{\mu s}(x_s,\la), \quad M_{sk\mu}(\la):=M_{sks\mu}(\la).  \eqno(11)
$$

We introduce the matrices $M_{s}(\la),\;s=\overline{1,p},$ as follows:
$$
M_{s}(\la)=[M_{sk\mu}(\la)]_{k,\mu=\overline{1,n}},\quad
M_{sk\mu}(\la):=\de_{k\mu}\; \mbox{for}\; k\ge \nu.
$$
The matrix $M_s(\la)$ is called the Weyl-type matrix with respect
to the boundary vertex $v_s$. The inverse problem is formulated
as follows. Fix $N=\overline{1,p}.$

\smallskip
{\bf Inverse problem 1.} Given $\{M_{s}(\la)\},\;
s=\overline{1,p}\setminus N$, construct $q$ on $T.$

\smallskip
We note that the notion of the Weyl-type matrices $M_s$ is a
generalization of the notion of the Weyl function (m-function) for
the classical Sturm-Liouville operator ([19, 27])  and is a
generalization of the notion of Weyl matrices introduced in [14,
15, 21, 22, 24] for higher-order differential operators on an
interval and on graphs. Thus, Inverse Problem 1 is a
generalization of the well-known inverse problems for differential
operators on an interval and on graphs.

\smallskip
We also note that in Inverse problem 1 we do not need to specify
all matrices $M_s(\lambda),$ $s=\overline{1,p}$; one of them can
be omitted. This last fact was first noticed in [8], where the
inverse problem was solved for the Sturm-Liouville operators on an
arbitrary tree.

\smallskip
In section 3 properties of the Weyl-type solutions  and the
Weyl-type matrices are studied. Section 4 is devoted to the
solution of auxiliary inverse problems of recovering the potential
on a fixed edge. In section 5 we study Inverse Problem 1. For this
inverse problem we provide a constructive procedure for the
solution and prove its uniqueness.\\

{\bf 3. Properties of spectral characteristics.} Fix
$s=\overline{1,p},\; k=\overline{1,n-1}.$ Substituting (10)-(11)
into matching conditions (7)-(8), we obtain a linear algebraic
system with respect to $M_{skj\mu}(\la).$ Solving this system by
Cramer's rule one gets
$M_{skj\mu}(\la)=\Delta_{skj\mu}(\la)/\Delta_{sk}(\la),$ where the
functions $\Delta_{skj\mu}(\la)$ and $\Delta_{sk}(\la)$ are entire
in $\la.$ Thus, the functions $M_{skj\mu}(\la)$ are meromorphic in
$\la,$ and consequently, the Weyl-type solutions and the Weyl-type
matrices are meromorphic in $\la.$ In particular,
$$
M_{sk\mu}(\la)=
\frac{\Delta_{sk\mu}(\la)}{\Delta_{sk}(\la)},\quad k\le\mu,           \eqno(12)
$$
where $\Delta_{sk\mu}(\la):=\Delta_{sks\mu}(\la).$ We note that
the function $\Delta_{sk}(\la)$ in (12) is the characteristic
function for the boundary value problem $L_{sk}$ for equation (1)
under the conditions
$$
y_{s}(x_s)=O(x_s^{\xi_{k+1,s}}),\; x_s\to 0,\qquad
y_{j}(x_j)=O(x_j^{\xi_{n-k+1,j}}),\; x_j\to 0,\; j=\overline{1,p},\;j\ne s,
$$
$$
U_{1\nu}(y_{1})=U_{j\nu}(y_{j}),\; j=\overline{2,p},\;\nu=\overline{0,k-1},
\qquad \sum_{j=1}^{p} U_{j\nu}(y_{j})=0,\;\nu=\overline{k,n-1}.
$$
Zeros of $\Delta_{sk}(\la)$ coincide with the eigenvalues of $L_{sk}.$
Denote
$$
\Omega_{kj}=\det[R_{l}^{\xi_{\mu j}}]_{l,\mu=\overline{1,k}},\;\Omega_{0j}:=1,
\quad \om_{kj}:=\frac{\Omega_{k-1,j}}{\Omega_{kj}},\; k=\overline{1,n}.
$$

{\bf Lemma 1. }{\it Fix $j=\overline{1,p},$ and fix a sector $S$
with property (2).

1) Let $k=\overline{1,n-1},$ and let $y_j(x_j,\la)$ be a solution
of equation (1) on the edge $e_j$ under the condition
$$
y_j(x_j,\la)=O(x_j^{\xi_{k+1,j}}),\quad x_j\to 0.                         \eqno(13)
$$
Then for $x_j\in(0,l_j],\;\nu=\overline{0,n-1},\;
\rho\in S,\; |\rho|\to\iy,$
$$
y_j^{(\nu)}(x_j,\la)=\sum_{\mu=k+1}^{n} A_{\mu j}(\rho)
(\rho R_{\mu})^{\nu}\exp(\rho R_{\mu} x_j)[1],                         \eqno(14)
$$
where the coefficients $A_{\mu j}(\rho)$ do not depend on $x_j$.
Here and below we assume that $\arg\rho=const,$ when $|\rho|\to\iy.$

2) Let $k=\overline{1,n},$ and let $y_j(x_j,\la)$ be a solution
of equation (1) on the edge $e_j$ under the condition
$$
y_j(x_j,\la)\sim c_{kj0}x_j^{\xi_{kj}},\; x_j\to 0.                       \eqno(15)
$$
Then for $x_j\in(0,l_j],\;\nu=\overline{0,n-1},\;\rho\in S,
\;|\rho|\to\iy,$
$$
y_j^{(\nu)}(x_j,\la)=\frac{\om_{kj}}{\rho^{\xi_{kj}}}
(\rho R_{k})^{\nu}\exp(\rho R_{k} x_j)[1]
+\sum_{\mu=k+1}^{n} B_{\mu j}(\rho)
(\rho R_{\mu})^{\nu}\exp(\rho R_{\mu} x_j)[1],                        \eqno(16)
$$
where the coefficients $B_{\mu j}(\rho)$ do not depend on $x_j$.}

\smallskip
{\bf Proof. } It follows from (13) that
$$
y_j(x_j,\la)=\sum_{\mu=k+1}^n a_{\mu j}(\la)S_{\mu j}(x_j,\la).       \eqno(17)
$$
Using the fundamental system of solutions
$\{E_{kj}(x_j,\rho)\}_{k=\overline{1,n}}$, one can write
$$
y_j(x_j,\la)=\sum_{m=1}^{n} A_{mj}(\rho)E_{mj}(x_j,\rho).             \eqno(18)
$$
By virtue of (3), we calculate
$$
y_j(x_j,\la)=\sum_{m=1}^{n} A_{mj}(\rho)\sum_{\mu=1}^n
b_{mj\mu}(\rho) S_{\mu j}(x_j,\la)=\sum_{\mu=1}^n S_{\mu j}(x_j,\la)
\sum_{m=1}^{n} A_{mj}(\rho) b_{mj\mu}(\rho).
$$
Comparing this relation with (17), we obtain
$$
\sum_{m=1}^{n} A_{mj}(\rho) b_{mj\mu}(\rho)=0,\quad \mu=\overline{1,k}.  \eqno(19)
$$
We consider (19) as a linear algebraic system with respect to
$A_{j}(\rho), A_{2j}(\rho),\ldots, A_{kj}(\rho).$
Solving this system by Cramer's rule and taking (4) into account we get
$$
A_{mj}(\rho)=\sum_{\mu=k+1}^{n}
(\al_{m\mu j}+O(\rho^{-1}))A_{\mu j}(\rho),\quad m=\overline{1,k},    \eqno(20)
$$
where $\al_{m\mu j}$ are constants. Substituting (20) into (18)
and using (2) we arrive at (14). Relations (16) are proved
analogously by using (15) instead of (13).
\B

\smallskip
Now we are going to study the asymptotic behavior of the Weyl-type solutions.

\smallskip
{\bf Lemma 2. }{\it Fix $s=\overline{1,p},\;k=\overline{1,n},$
and fix a sector $S$ with property (2). For $x_s\in(0,l_s),\;
\nu=\overline{0,n-1},$ the following asymptotic formula holds}
$$
\psi_{sks}^{(\nu)}(x_s,\la)=\frac{\om_{ks}}{\rho^{\xi_{ks}}}\,
(\rho R_{k})^{\nu}\exp(\rho R_{k} x_s)[1],\quad\rho\in S,\;|\rho|\to\iy. \eqno(21)
$$

{\bf Proof. } For $k=n$, (21) follows from Lemma 1.
Fix $s=\overline{1,p},\; k=\overline{1,n-1}.$
Using Lemma 1 and boundary conditions for $\Psi_{sk}$ we get the
following asymptotic formulae for $\rho\in S,\;|\rho|\to\iy$:
$$
\psi_{sks}^{(\nu)}(x_s,\la)=\frac{\om_{ks}}{\rho^{\xi_{ks}}}
\,(\rho R_{k})^{\nu}\exp(\rho R_{k} x_s)[1]+
\sum_{\mu=k+1}^{n} A_{\mu s}^{sk}(\rho)
(\rho R_{\mu})^{\nu}\exp(\rho R_{\mu} x_s)[1],\; x_s\in(0,l_s],        \eqno(22)
$$
$$
\psi_{skj}^{(\nu)}(x_j,\la)=\sum_{\mu=n-k+1}^{n}
A_{\mu j}^{sk}(\rho)(\rho R_{\mu})^{\nu}\exp(\rho R_{\mu} x_j)[1],
\;j=\overline{1,p}\setminus s,\;  x_j\in(0,l_j].                      \eqno(23)
$$
Substituting (22)-(23) into matching conditions (7)-(8) for
$\Psi_{sk}$, we obtain the linear algebraic system with respect to
$A_{\mu j}^{sk}(\rho).$ Solving this system by Cramer's rule,
we obtain in particular,
$$
A_{\mu j}^{sk}(\rho)=O(\rho^{-\xi_{ks}}\exp(\rho(R_{k}-R_{\mu})l_j)). \eqno(24)
$$
Substituting (24) into (22) we arrive at (21).
\B

\smallskip
It follows from the proof of Lemma 2 that one can also get the
asymptotics for $\psi_{skj}^{(\nu)}(x_j,\la),$ $j\ne s$; but for
our purposes only (21) is needed.\\

{\bf 4. Auxiliary inverse problems. } In this section we consider
auxiliary inverse problems of recovering differential operator on
each fixed edge. Fix $s=\overline{1,p},$ and consider the
following inverse problem on the edge $e_s$.

\smallskip
{\bf IP(s). } Given the Weyl-type matrix $M_s$,
construct the potential $q_{s}$ on the edge $e_s$.

\smallskip
In this inverse problem we construct the potential only
on the edge $e_s$, but the Weyl-type matrix $M_s$ brings a global
information from the whole graph. In other words, this problem
is not a local inverse problem related only to the edge $e_s$.

\smallskip
Let us formulate the uniqueness theorem for the solution of the
inverse problem $IP(s)$. For this purpose together with $q$ we
consider a potential $\tilde q.$ Everywhere below if a symbol
$\al$ denotes an object related to $q,$ then $\tilde\al$ will
denote the analogous object related to $\tilde q.$

\smallskip
{\bf Theorem 1. }{\it Fix $s=\overline{1,p}.$ If $M_s=\tilde M_s,$
then $q_{s}=\tilde q_{s}.$ Thus, the specification of the Weyl-type
matrix $M_s$ uniquely determines the potential $q_s$ on the edge $e_s$.}

\smallskip
We omit the proof since it is similar to that in [22, Ch.2].
Moreover, using the method of spectral mappings and Lemma 2, one
can get a constructive procedure for the solution of the inverse
problem $IP(s).$ It can be obtained by the same arguments as for
$n$-th order differential operators on a finite interval (see [22,
Ch.2] for details). Note that like in [22], the nonlinear inverse
problem $IP(s)$ is reduced to the solution of a linear equation in
the corresponding Banach space of sequences. The unique
solvability of this linear equation is proved by the same
arguments as in [22].

\medskip
Fix $j=\overline{1,p}.$ Now we define an auxiliary Weyl-type matrix
with respect to the internal vertex $v_0$ and the edge $e_j$.
Let $\vv_{kj}(x_j,\la),$ $k=\overline{1,n},$ be solutions of
equation (1) on the edge $e_j$ under the conditions
$$
\vv_{kj}^{(\nu-1)}(l_j,\la)=\de_{k\nu},\;\nu=\overline{1,k},
\qquad \vv_{kj}(x_j,\la)=O(x_j^{\xi_{n-k+1,j}}),\;x_j\to 0.
$$
We introduce the matrix
$m_j(\la)=[m_{jk\nu}(\la)]_{k,\nu=\overline{1,n}},$
where $m_{jk\nu}(\la):=\vv^{(\nu-1)}_{kj}(l_j,\la).$
Clearly, $m_{jk\nu}(\la)=\de_{k\nu}$ for $k\ge\nu,$ and
$\det m_j(\la)\equiv 1.$
The matrix $m_j(\la)$ is called the Weyl-type matrix with
respect to the internal vertex $v_0$ and the edge $e_j$.
Consider the following inverse problem on the edge $e_j$.

\smallskip
{\bf IP[j]. } Given the Weyl-type matrix $m_j$, construct
the potential $q_{j}$ on the edge $e_j$.

\smallskip
This inverse problem is the classical one, since it is the inverse
problem of recovering a higher-order differential equation on a
finite interval from its Weyl-type matrix. This inverse problem
has been solved in [22], where the uniqueness theorem is proved.
Moreover, in [22] an algorithm for the solution of the inverse
problem $IP[j]$ is given, and necessary and sufficient conditions
for the solvability of this inverse
problem are provided.\\

{\bf 5. Solution of Inverse Problem 1. } In this section we obtain
a constructive procedure for the solution of Inverse problem 1 and
prove its uniqueness. First we prove an auxiliary assertion.

\smallskip
{\bf Lemma 3. }{\it Fix $j=\overline{1,p}.$
Then for each fixed $s=\overline{1,p}\setminus j,$}
$$
m_{j1\nu}(\la)=\di\frac{\psi_{s1j}^{(\nu-1)}(l_j,\la)}{\psi_{s1j}(l_j,\la)},
\quad \nu=\overline{2,n},                                                  \eqno(25)
$$
$$
m_{jk\nu}(\la)=\di\frac{\det[\psi_{s\mu j}(l_j,\la),\ldots,
\psi_{s\mu j}^{(k-2)}(l_j,\la),
\psi_{s\mu j}^{(\nu-1)}(l_j,\la)]_{\mu=\overline{1,k}}}
{\det[\psi_{s\mu j}^{(\xi-1)}(l_j,\la)]_{\xi,\mu=\overline{1,k}}}\,,
\;2\le k<\nu\le n.                                                        \eqno(26)
$$

{\bf Proof. } Denote
$$
w_{js}(x_j,\la):=\di\frac{\psi_{s1j}(x_j,\la)}{\psi_{s1j}(l_j,\la)}\,.
$$
The function $w_{js}(x_j,\la)$ is a solution of equation (1) on the edge
$e_j$, and $w_{js}(l_j,\la)=1.$ Moreover, by virtue of the boundary
conditions on $\Psi_{s1}$, one has $w_{js}(x_j,\la)=O(x_j^{\xi_{nj}}),$
$x_j\to 0.$ Hence, $w_{js}(x_j,\la)\equiv\vv_{1j}(x_j,\la),$
i.e.
$$
\vv_{1j}(x_j,\la)=\di\frac{\psi_{s1j}(x_j,\la)}{\psi_{s1j}(l_j,\la)}.       \eqno(27)
$$
Similarly, we calculate
$$
\vv_{kj}(x_j,\la)=\di\frac{\det[\psi_{s\mu j}(l_j,\la),\ldots,
\psi_{s\mu j}^{(k-2)}(l_j,\la),\psi_{s\mu j}(x_j,\la)]_{\mu=\overline{1,k}}}
{\det[\psi_{s\mu j}^{(\xi-1)}(l_j,\la)]_{\xi,\mu=\overline{1,k}}}\,,
\quad k=\overline{2,n-1}.                                                 \eqno(28)
$$
Since $m_{jk\nu}(\la)=\vv_{kj}^{(\nu-1)}(l_j,\la),$ it follows
from (27) that (25) holds. Similarly, (26) follows from (28).
\B

\medskip
Now we are going to obtain a constructive procedure for the
solution of Inverse problem 1. Our plan is the following.

{\it Step 1. } Let the Weyl-type matrices $\{M_{s}(\la)\},\;
s=\overline{1,p}\setminus N$, be given. Solving the inverse problem
$IP(s)$ for each fixed $s=\overline{1,p}\setminus N,$ we find the
potentials $q_{s}$ on the edges $e_s$, $s=\overline{1,p}\setminus N$.

{\it Step 2. } Using the knowledge of the potential on the edges
$e_s$, $s=\overline{1,p}\setminus N$, we construct the Weyl-type
matrix $m_{N}(\la)$.

{\it Step 3. } Solving the inverse problem $IP[N]$, we find
the potential $q_{N}$ on $e_{N}$.

\smallskip
Steps 1 and 3 have been already studied in Section 4.
It remains to fulfil Step 2.

Suppose that Step 1 is already made, and we found the potentials
$q_{s}$, $s=\overline{1,p}\setminus N$, on the edges $e_s$,
$s=\overline{1,p}\setminus N$. Then we calculate the functions
$S_{kj}(x_j,\la),$ $j=\overline{1,p}\setminus N,$
$k=\overline{1,n}.$

Fix $s=\overline{1,p}\setminus N.$
All calculations below will be made for this fixed $s.$

Our goal now is to construct the Weyl-type matrix $m_{N}(\la).$
For this purpose we will use Lemma 3. According to (25)-(26), in
order to construct $m_{N}(\la)$ we have to calculate the functions
$$
\psi_{skN}^{(\nu)}(l_{N},\la),\quad k=\overline{1,n-1},\;\nu=\overline{0,n-1}. \eqno(29)
$$
We will find the functions (29) by the following steps.

\smallskip
1) Using (11) we construct the functions
$$
\psi_{sks}^{(\nu)}(l_s,\la),\;k=\overline{1,n-1},\;\nu=\overline{0,n-1},     \eqno(30)
$$
by the formula
$$
\psi_{sks}^{(\nu)}(l_s,\la)=S_{ks}^{(\nu)}(l_s,\la)+
\sum_{\mu=k+1}^{n} M_{sk\mu}(\la)S_{\mu s}^{(\nu)}(l_s,\la).              \eqno(31)
$$

\smallskip
2) Using the matching conditions (7) on $\Psi_{sk}$, we get, in particular,
$$
U_{j\nu}(\psi_{skj})=U_{s\nu}(\psi_{sks}),
\quad 0\le\nu<k\le n-1,\quad j=\overline{1,p}\setminus s.                   \eqno(32)
$$
Since the functions (30) were already calculated, it follows that the
right-hand sides in (32) are known. For each fixed $k=\overline{1,n-1},$
we successively use (32) for $\nu=0,1,\ldots, k-1,$ and calculate
recurrently the functions
$$
\psi_{skj}^{(\nu)}(l_j,\la),\;k=\overline{1,n-1},\;
\nu=\overline{0,k-1},\; j=\overline{1,p}\setminus s.                       \eqno(33)
$$
In particular we found the functions (29) for $\nu=\overline{0,k-1}.$

\smallskip
3) It follows from (10) that
$$
\sum_{\mu=n-k+1}^{n} M_{skj\mu}(\la)S_{\mu j}^{(\nu)}(l_j,\la)
=\psi_{skj}^{(\nu)}(l_j,\la),\;k=\overline{1,n-1},
\;j=\overline{1,p}\setminus s,\;\nu=\overline{0,n-1}.                     \eqno(34)
$$
Fix $k=\overline{1,n-1},\; j=\overline{1,p},\;j\ne s,\; j\ne N,$ and
consider a part of relations (34), namely, for $\nu=\overline{0,k-1}.$
For this choice of the parameters, the right-hand sides in (34) are
known, since the functions (33) are known. Relations (34) for
$\nu=\overline{0,k-1},$ form a linear algebraic system $\sigma_{skj}$
with respect to the coefficients $M_{skj\mu}(\la),\;
\mu=\overline{n-k+1,n}.$ Solving this system by Cramer's rule, we find
this functions. Substituting them into (34), we calculate the functions
$$
\psi_{skj}^{(\nu)}(l_j,\la),\quad k=\overline{1,n-1},
\;j=\overline{1,p}\setminus N,\; \nu=\overline{0,n-1}.                  \eqno(35)
$$
Note that for $j=s$ these functions were found earlier (see(31)).

\smallskip
4) Let us now use the generalized Kirchhoff's conditions (8) for $\Psi_{sk}$.
Since the functions (35) are known, one can construct by (8) the functions
(29) for $k=\overline{1,n-1},\;\nu=\overline{k,n-1}.$ Thus, the functions
(29) are known for $k=\overline{1,n-1},\; \nu=\overline{0,n-1}.$

\smallskip
Since the functions (29) are known, we construct the Weyl-type matrix
$m_{N}(\la)$ via (25)-(26) for $j=N.$ Thus, we have obtained the solution
of Inverse problem 1 and proved its uniqueness, i.e. the following
assertion holds.

\medskip
{\bf Theorem 2. }{\it The specification of the Weyl-type matrices $M_s(\la),$
$s=\overline{1,p}\setminus N$, uniquely determines the potential $q$ on $T.$
The solution of Inverse problem 1 can be obtained by the following algorithm.}

\medskip
{\bf Algorithm 1. }{\it Given the Weyl-type matrices
$M_s(\la),$ $s=\overline{1,p}\setminus N$.

1) Find the potentials $q_{s}$, $s=\overline{1,p}\setminus N$, by solving
the inverse problem $IP(s)$ for each fixed $s=\overline{1,p}\setminus N$.

2) Calculate $S_{kj}^{(\nu)}(l_j,\la),\; j=\overline{1,p}\setminus N,$
$k=\overline{1,n},\;\nu=\overline{0,n-1}.$

3) Fix $s=\overline{1,p}\setminus N$. All calculations below will be made
for this fixed $s.$ Construct the functions (30) via (31).

4) Calculate the functions (33) using (32).

5) Find the functions $M_{skj\mu}(\la),$ by solving the linear algebraic
systems $\sigma_{skj}$.

6) Construct the functions (35) using (34).

7) Find the functions (29) using (33), (35) and (8).

8) Calculate the Weyl-type matrix $m_{N}(\la)$ via (25)-(26) for
$j=N$.

9) Construct the potential $q_{N}$ on the edge $e_{N}$ by solving
the inverse problem $IP[N]$.}

\bigskip
\noindent {\bf Acknowledgment.}  This work was supported by Grant 1.1436.2014K of the Russian Ministry of Education and Science and by Grant 13-01-00134 of Russian Foundation for Basic Research.

\begin{center}
{\bf REFERENCES}
\end{center}
\begin{enumerate}
\item[{[1]}] E. Montrol, {\it Quantum theory on a network}, J. Math. Phys. 11, no.2 (1970), 635-648.
\item[{[2]}] J. Langese, G. Leugering and J. Schmidt, {\it Modelling, analysis and control
     of dynamic elastic multi-link structures}, Birkh\"auser, Boston, 1994.
\item[{[3]}] T. Kottos and U. Smilansky, {\it Quantum chaos on graphs}, Phys. Rev. Lett. 79 (1997), 4794-4797.
\item[{[4]}] P. Kuchment, {\it Quantum graphs. Some basic structures}, Waves Random Media 14 (2004), S107-S128.
\item[{[5]}] Yu. Pokornyi and A. Borovskikh, {\it Differential equations on  networks (geometric graphs)}, J. Math. Sci. (N.Y.) 119, no.6 (2004), 691-718.
\item[{[6]}] Yu. Pokornyi and V. Pryadiev, {\it The qualitative Sturm-Liouville theory on spatial networks}, J. Math. Sci. (N.Y.) 119 (2004), no.6, 788-835.
\item[{[7]}] M.I. Belishev, {\it Boundary spectral inverse problem on a class of graphs (trees) by the BC method}, Inverse Problems 20 (2004), 647-672.
\item[{[8]}] V.A. Yurko, {\it Inverse spectral problems for Sturm-Liouville operators on graphs}, Inverse Problems 21 (2005), 1075-1086.
\item[{[9]}] B.M. Brown and R. Weikard, {\it A Borg-Levinson theorem for trees}, Proc. R. Soc. Lond. Ser. A Math. Phys. Eng. Sci. 461, no.2062 (2005), 3231-3243.
\item[{[10]}] V.A. Yurko, {\it Inverse problems for Sturm-Liouville operators on bush-type graphs}, Inverse Problems 25, no.10 (2009), 105008, 14pp.
\item[{[11]}] V.A. Yurko, {\it An inverse problem for Sturm-Liouville operators on A-graphs}, Applied Math. Lett. 23, no.8 (2010), 875-879.
\item[{[12]}] V.A. Yurko, {\it Inverse spectral problems for differential operators on arbitrary compact graphs}, Journal of Inverse and Ill-Posed Proplems 18, no.3 (2010), 245-261.
\item[{[13]}] C.-F. Yang, {\it Inverse spectral problems for Sturm-Liouville operators on a d-star graph}, J. Math. Anal. Appl. 365 (2010), 742-749.
\item[{[14]}] V.A. Yurko, {\it An inverse problem for higher-order differential operators on star-type graphs}, Inverse Problems 23, no.3 (2007), 893-903.
\item[{[15]}] V.A. Yurko, {\it Inverse problems for differential of any order on trees}, Matemat. Zametki 83, no.1 (2008), 139-152; English  transl. in Math. Notes 83, no.1 (2008), 125-137.
\item[{[16]}] V.A. Marchenko, {\it Sturm-Liouville operators and their applications}, ''Naukova Dumka",  Kiev, 1977;  English  transl., Birkh\"auser, 1986.
\item[{[17]}] B.M. Levitan, {\it Inverse Sturm-Liouville problems}, Nauka, Moscow, 1984; English transl., VNU Sci.Press, Utrecht, 1987.
\item[{[18]}] K. Chadan, D. Colton, L. Paivarinta and W. Rundell, {\it An introduction to inverse scattering and inverse spectral problems}, SIAM Monographs on Mathematical Modeling and Computation. SIAM, Philadelphia, PA, 1997.
\item[{[19]}] G. Freiling and V.A. Yurko, {\it Inverse Sturm-Liouville Problems and their Applications}, NOVA Science Publishers, New York, 2001.
\item[{[20]}] R. Beals, P. Deift and C. Tomei, {\it Direct and Inverse Scattering on the Line}, Math. Surveys and Monographs, v.28. Amer. Math. Soc. Providence: RI, 1988.
\item[{[21]}] V.A. Yurko, {\it Inverse Spectral Problems for Differential Operators and their Applications}, Gordon and Breach, Amsterdam, 2000.
\item[{[22]}] V.A. Yurko, {\it Method of Spectral Mappings in the Inverse Problem Theory}, Inverse and Ill-posed Problems Series. VSP, Utrecht, 2002.
\item[{[23]}] V.A. Yurko, {\it Inverse  problem for differential equations with a singularity}, Differ. Uravneniya, 28, no.8 (1992), 1355-1362 (Russian); English transl. in Diff. Equations, 28 (1992), 1100-1107.
\item[{[24]}] V.A. Yurko, {\it On higher-order differential operators with a singularity}, Matem. Sbornik, 186, no.6 (1995), 133-160 (Russian); English transl. in Sbornik; Mathematics 186, no.6 (1995), 901-928.
\item[{[25]}] V.A. Yurko, {\it Inverse spectral problems for higher-order differential operators with a singularity}, Journal of Inverse and Ill-Posed Problems, 10, no.4 (2002), 413-425.
\item[{[26]}] V.A. Yurko, {\it Higher-order differential equations having a singularity in an interior point}, Results in Mathematics, 42, no.1-2 (2002), 177-191.
\item[{[27]}] B.M. Levitan and I.S. Sargsyan, {\it Introduction to Spectral Theory}, AMS Transl. of Math. Monogr. 39, Providence, 1975.
\end{enumerate}

\medskip
\begin{tabular}{ll}
Name:             &   Yurko, Vjacheslav  \\
Place of work:    &   Department of Mathematics, Saratov State University \\
{}                &   Astrakhanskaya 83, Saratov 410012, Russia \\
Present Position: &   Professor, Head of the Faculty of Mathematical Physics \\
Telephone:        &   (8452) 275526           \\
E-mail:           &   yurkova@info.sgu.ru
\end{tabular}

\end{document}